\documentclass[12pt]{article}
\usepackage{epic,latexsym,amssymb}
\usepackage{color}
\usepackage{tikz}

\usepackage{geometry,graphicx,verbatim,amsmath}
\usepackage{tikz,ifthen}
\usetikzlibrary{calc}

\newtheorem{thm}{Theorem}[section]
\newtheorem{lem}[thm]{Lemma}
\newtheorem{cor}[thm]{Corollary}

\newtheorem{prop}[thm]{Proposition}

\newcommand{\mdim}{{\rm mdim}}

\newcommand{\CV}{{\rm CV}}

\newcommand{\proof}{\noindent{\bf Proof.\ }}
\newcommand{\qed}{\hfill $\square$ \bigskip}

\begin{document}

\title{Graphs whose mixed metric dimension is equal to their order}

\author{ Ali Ghalavand$^{a}$ \and Sandi Klav\v{z}ar$^{b,c,d}$ \and Mostafa Tavakoli$^{a,}$\footnote{Corresponding author}} 

\date{}

\maketitle
\vspace{-0.8 cm}
\begin{center}
$^a$ Department of Applied Mathematics, Faculty of Mathematical Sciences,\\
Ferdowsi University of Mashhad, P.O.\ Box 1159, Mashhad 91775, Iran\\
{\tt m$\_$tavakoli@um.ac.ir\\ 
alighalavand@grad.kashanu.ac.ir}\\
\medskip

$^b$ Faculty of Mathematics and Physics, University of Ljubljana, Slovenia\\
{\tt sandi.klavzar@fmf.uni-lj.si}\\
\medskip

$^c$ Faculty of Natural Sciences and Mathematics, University of Maribor, Slovenia\\
\medskip

$^d$ Institute of Mathematics, Physics and Mechanics, Ljubljana, Slovenia\\
\end{center}

\begin{abstract}
The mixed metric dimension ${\rm mdim}(G)$ of a graph $G$ is the cardinality of a smallest set of vertices that (metrically) resolves each pair of elements from $V(G)\cup E(G)$. We say that $G$ is a  max-mdim graph if ${\rm mdim}(G) = n(G)$. It is proved that a max-mdim graph $G$ with $n(G)\ge 7$ contains a vertex of degree at least $5$. Using the strong product of graphs and amalgamations large families of max-mdim graphs are constructed. The mixed metric dimension of graphs with at least one universal vertex is determined. The mixed metric dimension of graphs $G$ with cut vertices is bounded from the above and the mixed metric dimension of block graphs computed. 
\end{abstract}


\newpage

\section{Introduction}

The metric dimension is an extremely prolific and at the same time interesting area of graph theory, for several reasons. The main reason is certainly that the theory is extremely useful in other areas of science, for instance in computer science, chemistry, social networks, and biology, see respective papers~\cite{Melter-1984, Johnson-1993, Trujillo-Rasua-2016, till-2019}. For more information on the metric dimension and its applications see the recent survey \cite{till-2022+}. On the other hand, various applications also give rise to certain modifications of the basic concept, which leads to further intensive research to obtain additional insight into the classical topic and between variants. For more information on this point of view of the metric dimension see the other recent survey \cite{dorota-2022+}.  We also refer to a recent application of the local metric dimension to delivery services from~\cite{KT-2021}. 

A very interesting version of the metric dimension was introduced in 2017 by Kelenc, Kuziak, Taranenko, and Yero~\cite{KKTY-2017}, namely mixed metric dimension, as follows. Let $G=(V(G), E(G))$ be a graph. Then two elements $x,y\in V(G)\cup E(G)$ are \emph{resolved} by a vertex $v\in V(G)$ if $d_G(x,v)\ne d_G(y,v)$, where $d_G$ stands either for the shortest-path distance between vertices, or the distance between an edge and a vertex. The latter distance is, for an edge $x=ww'$ and a vertex $v$, defined by $d_G(x,v)=\min\{d_G(w,v), d_G(w',v)\}$. A set of vertices $W\subset V(G)$ is a \emph{mixed resolving set} for $G$ if any two elements (vertices or edges) $x,y\in V(G)\cup E(G)$ are resolved by a vertex of $W$.
We note that $V(G)$ is always a mixed resolving set for $G$.
 A mixed resolving set of the smallest cardinality is  a mixed metric basis, its cardinality is the mixed metric dimension $\mdim(G)$. After the seminal paper, the mixed metric dimension was investigated in many papers, cf.~\cite{ghalanand-2022+, milivojevic-2021, MKSM-2021, qu-2022, raza-2020, sharma-2022, SS-2021}.

Let $G$ be a graph and $x\in V(G)$. Then a neighbor $y$ of $x$ is a {\em maximal neighbor of} $x$ if $y$ is adjacent to all neighbors of $x$. Denoting the order of $G$ by $n(G)$, we recall the following result which is the main source of inspiration for this article. 

\begin{thm}\label{thm:motivation} {\rm \cite[Theorem 3.8]{KKTY-2017} }
If $G$ is a graph, then $\mdim(G) = n(G)$ if and only if every vertex of $G$ has a maximal neighbour.
\end{thm}

Let us say that a graph $G$ with $\mdim(G) = n(G)$ is a {\em max-mdim graph}. Theorem~\ref{thm:motivation} thus characterizes max-mdim graphs. The main purpose of this article is to take a closer look at this class of graphs. We proceed as follows. In the rest of the introduction some further definitions are listed and a result is stated to be used later on. In the next section we prove that a max-mdim graph $G$ with $n(G)\ge 7$ contains a vertex of degree at least $5$. This implies that if $G$ is a chemical graph, then $\mdim(G)\le n(G)-1$. Afterwards we apply the strong product and amalgamations to construct large families of max-mdim graphs. In particular, the strong products $P_n\boxtimes K_2$ are max-mdim graphs with $\Delta = 5$ where $\Delta$ is 
the largest number of neighbors of a vertex in $P_n\boxtimes K_2$.
We also determine the mixed metric dimension for graphs with universal vertices. In the concluding section we consider the mixed metric dimension of graphs $G$ with cut vertices and prove an upper bound on their mixed metric dimension.  As a consequence we determine the  mixed metric dimension of block graphs. 

In this paper, we consider finite, simple and connected graphs. Let $G$ be a graph. The degree of $v\in V(G)$ will be denoted by $\deg_G(v)$. The (open) neighborhood of $v$ will be denoted by $N_G(v)$. Then $\deg_G(v) = |N_G(v)|$. A pendant vertex of $G$ is a vertex with degree one. A vertex of degree $n(G)-1$ is a {\em universal vertex}. The minimum and the maximum degree of $G$ are respectively denoted by $\delta(G)$ and $\Delta(G)$. The number of cut vertices of a $G$ is denoted by $\zeta(G)$ and the set of all cut vertices by $\CV(G)$, so that $|\CV(G)| = \zeta(G)$.
A {\em block} of a graph is a nonseparable maximal subgraph of the graph. A graph is $2$-{\em connected} if it has no cut vertices.
 Note that if $G$ is $2$-connected, then $\zeta(G) = 0$. $G$ is a {\em block graph} if each block of $G$ is complete. To conclude this article introduction, we state the following result which implicitly follows from~\cite[Theorem 3.8]{KKTY-2017}.

\begin{lem}\label{lem:max-neigbor}
If $W$ is a mixed resolving set for a graph $G$, and $v\in V(G)$ has a maximal neighbor, then $v\in W$.
\end{lem}

\section{Classes of max-mdim graphs and a maximum degree bound}
\label{sec:max-degree}

In this section we prove that max-mdim graphs necessarily contain a vertex of degree at least $5$ as soon as their order is at least $7$. Then we use the strong product and amalgamations to construct large families of max-mdim graphs. In particular, the strong products $P_n\boxtimes K_2$ are max-mdim graphs with $\Delta =5$. We also determine the mixed metric dimension for graphs with universal vertices. 

\begin{figure}[ht!]
\begin{center}
\begin{tikzpicture}[scale=0.8,style=thick]
\tikzstyle{every node}=[draw=none,fill=none]
\def\vr{3pt} 
\path (0,0) coordinate (v1);
\path (2,0) coordinate (v2);
\path (4,0) coordinate (v3);
\path (0,2) coordinate (v4);
\path (2,2) coordinate (v5);
\path (4,2) coordinate (v6);
\draw (v1) -- (v2) -- (v3) -- (v6) -- (v5) -- (v4)  -- (v1);
\draw (v1) -- (v5) -- (v3);
\draw (v4) -- (v2) -- (v6);
\draw (v5) -- (v2);
%
\draw (v1)  [fill=black] circle (\vr);
\draw (v2)  [fill=black] circle (\vr);
\draw (v3)  [fill=black] circle (\vr);
\draw (v4)  [fill=black] circle (\vr);
\draw (v5)  [fill=black] circle (\vr);
\draw (v6)  [fill=black] circle (\vr);
\end{tikzpicture}
\end{center}
\caption{A max-mdim graph}
\label{fig:P3-strong-P2}
\end{figure}
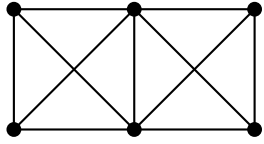

Since every vertex of a complete graph has a maximal neighbour,  by Theorem \ref{thm:motivation},  complete graphs are max-mdim graphs.
 Moreover, if $G$ is obtained from a complete graph by removing a matching, then $G$ is also a max-mdim graph provided that it contains at least two universal vertices. In particular, $K_4-e$ is  a max-mdim graph. Another small example is shown in Fig~\ref{fig:P3-strong-P2}. 

In our first theorem we prove that a max-mdim graph contains a vertex of degree at least $5$ as soon as it is not very small. 

\begin{thm}\label{thm:Delta}
If $G$ is a (connected) max-mdim graph with $n(G)\ge 7$, then $\Delta(G) \ge 5$.
\end{thm}

\proof
Let $G$ be a max-mdim graph with $n(G)\ge 7$ and set $\Delta = \Delta(G)$ for the rest of the proof. 

Suppose $\Delta = 2$. If $\delta(G) = 1$, the support vertex of a pendant vertex does not admit a maximal neighbor. Otherwise $G$ is a cycle which is not a max-mdim graph. Suppose $\Delta = 3$. Let $x$ be a vertex of $G$ with $\deg(x) = 3$ and let $N_G(x) = \{x_1, x_2, x_3\}$. Without loss generality assume that $x_1$ is a maximal neighbor of $x$, so that $x_1x_2, x_1x_3\in E(G)$ and $\deg_G(x_1) = 3$. As $n(G)\ge 7$ and $G$ is connected,  $x_2$ or $x_3$ is also of degree  $3$. Assume without loss of generality that $N_G(x_2) = \{x, x_1, x_2'\}$. Now, no matter which neighbor of $x_2$ is its maximal neighbor, we get $\deg_G(x) \ge 4$ or $\deg_G(x_1)\ge 4$, which is not possible. 

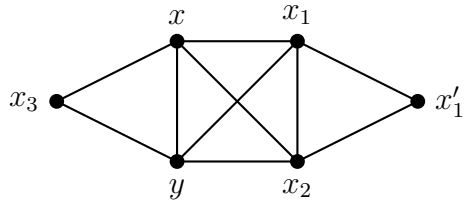
\begin{figure}[ht!]
\begin{center}
\begin{tikzpicture}[scale=0.8,style=thick]
\tikzstyle{every node}=[draw=none,fill=none]
\def\vr{3pt} 
\path (0,1) coordinate (v1);
\path (2,0) coordinate (v2);
\path (2,2) coordinate (v3);
\path (4,0) coordinate (v4);
\path (4,2) coordinate (v5);
\path (6,1) coordinate (v6);
\draw (v1) -- (v2) -- (v4) -- (v6) -- (v5) -- (v3)  -- (v1);
\draw (v3) -- (v2) -- (v5) -- (v4) -- (v3);
%
\draw (v1)  [fill=black] circle (\vr);
\draw (v2)  [fill=black] circle (\vr);
\draw (v3)  [fill=black] circle (\vr);
\draw (v4)  [fill=black] circle (\vr);
\draw (v5)  [fill=black] circle (\vr);
\draw (v6)  [fill=black] circle (\vr);
\draw[below] (v2)++(0.0,-0.1) node {$y$};
\draw[below] (v4)++(0.0,-0.1) node {$x_2$};
\draw[above] (v3)++(0.0,+0.1) node {$x$};
\draw[above] (v5)++(0.0,+0.1) node {$x_1$};
\draw[left] (v1)++(-0.1,0.0) node {$x_3$};
\draw[right] (v6)++(0.1,0.0) node {$x_1'$};
\end{tikzpicture}
\end{center}
\caption{The graph $G_6$}
\label{fig:G6}
\end{figure}

Suppose $\Delta = 4$. Let $x$ be a vertex of $G$ with $\deg(x) = 4$ and let $y$ be its maximal neighbor. Then $N_G(x) = N_G(y)$, say $N_G(x) = N_G(y) = \{x_1, x_2, x_3\}$.  As $n(G)\ge 7$, there is another vertex of $G$, without loss of generality assume it is adjacent to $x_1$, denote it by $x_1'$. Consider now a maximal neighbor of $x_1$. It cannot be $x$ or $y$ because then $x$ or $y$ would be adjacent to $x_1'$ and hence $x$ or $y$ would be of degree at least $5$. 
For the same reason, a maximal neighbor of $x_1$ cannot be $x_1'$. If there were another neighbor $x_1''$ of $x_1$, it also cannot be a maximal neighbor of $x_1$. So $x_1$ must have a maximal neighbor among the already introduced vertices and thus the fourth neighbor of $x_1$ is from $\{x_2, x_3\}$. If $x_1x_3\in E(G)$, then $x_3$ is the only candidate for a maximal neighbor of $x_1$ and therefore $x_1'x_3\in E(G)$, while if $x_1x_2\in E(G)$, then we must have $x_1'x_2\in E(G)$. But in both cases we get isomorphic graphs, see Fig.~\ref{fig:G6}, where the labeling presented is with respect to the second case.  

Hence, if $\Delta = 4$, then $G$ necessary contains the graph $G_6$ from Fig.~\ref{fig:G6} as an induced subgraph. Since $n(G)\ge 7$, there is another vertex, say $x_1''$, and we may assume without loss of generality that $x_1'x_1''\in E(G)$. We now infer that none of the vertices $x_1$, $x_2$, $x_1''$, or a possible fourth new neighbor of $x_1'$ can be a maximal neighbor of $x_1'$. So the only possibility is that $x_1'x_3\in E(G)$ so that $x_3$ would be  a maximal neighbor of $x_1'$. But then $x_3x_1, x_2x_1\in E(G)$ which means that $x_1$ would be of degree at least $5$.  
\qed

Note that the proof of Theorem~\ref{thm:Delta} implies that the graph $G_6$ from Fig.~\ref{fig:G6} is the unique max-mdim graph with $n(G) = 6$ and $\Delta(G) = 4$.  

In view of the applicability of the mixed metric dimension in chemistry~\cite{sharma-2022}, we recall that a graph $G$ is called a {\em chemical graph} if $\Delta(G) \le 4$. Theorem~\ref{thm:Delta} cleary has the following application. 

\begin{cor}
If $G$ is a chemical graph with $n(G) \ge 7$, then $\mdim(G) \le n(G) - 1$. 
\end{cor}

The graphs in Figs.~\ref{fig:P3-strong-P2} and~\ref{fig:G6} motivate us to recall the definition of the {\em strong product} $G\boxtimes H$ of graphs $G$ and $H$. Its vertex set is $V( G\boxtimes H) = V(G)\times V(H)$, and vertices $(g_1,h_1)$ and $(g_2,h_2)$ are adjacent if $h_1=h_2$ and $g_1$ is adjacent to $g_2$, or $g_1=g_2$ and $h_1$ is adjacent to $h_2$, or $g_1$ is adjacent to $g_2$ and $h_1$ is adjacent to $h_2$. A standard reference for the strong product is the book~\cite{hik-2011}. The metric dimension of strong products was investigated in~\cite{rodr-2015}, and the local metric dimension of strong products in~\cite{barr-2016}. We now use this graph operation to significantly increase the variety of max-mdim graphs. 

\begin{prop}
\label{prop:strong}
If $G$ is a graph, then $G\boxtimes K_2$ is a max-mdim graph. 
\end{prop}

\proof
Let $n = n(G)$, let $V(G) = \{v_1, \ldots, v_n\}$ and $V(K_2) = \{0,1\}$. Then by the definition of the strong product, $(v_i,1)$ is a maximal neighbor of $(v_i,0)$, and $(v_i,0)$ is a maximal neighbor of $(v_i,1)$ for each $i\in [n]$. Therefore, $G\boxtimes K_2$ is a max-mdim graph by Theorem~\ref{thm:motivation}. 
\qed

The special case $P_n\boxtimes K_2$ of Proposition~\ref{prop:strong} gives an infinite family of max-mdim graphs $G$ with $\Delta(G) = 5$. Hence Theorem~\ref{thm:Delta} cannot be improved in general. 

Another source for max-mdim graphs is the following construction. Let $G$ and $H$ be disjoint graphs, $e_G \in E(G)$ and $e_H\in E(H)$. Then the graph $A(G,e_G;H,e_H)$ is obtained from the disjoint union of $G$ and $H$ by identifying the edges $e_G$ and $e_H$. (``A" stands here for an amalgamation.) Actually, this identification can be done in two ways, but for our purposes any of these will do it.  

\begin{prop}
\label{prop:amalg}
If $G$ and $H$ are max-mdim graphs, and $e_G$ and $e_H$ are edges whose endpoints are maximal neighbors for each other, then $A(G,e_G;H,e_H)$ is a max-mdim graph. 
\end{prop}

\proof
Set $A = A(G,e_G;H,e_H)$. Let $e_G = gg'$ and $e_H = hh'$. If $x\in V(G)\setminus \{g,g'\}$, then its maximal neighbor in $G$ is also a maximal neighbor of $x$ in $A$. Similarly, a maximal neighbor of $y\in V(H)\setminus \{h,h'\}$ is a maximal neighbor of $y$ in $A$. Finally, $g=h$ is a maximal neighbor of $g' = h'$ in $A$, and $g'=h'$ is a maximal neighbor of $g = h$.
\qed

Using Proposition~\ref{prop:amalg}, we can state the following result. 
 
\begin{thm}
If $n \geq t \geq 5$, then there exists a max-mdim graph $G$ with $n(G) = n$ and $\Delta(G)=t$. 
 \end{thm}

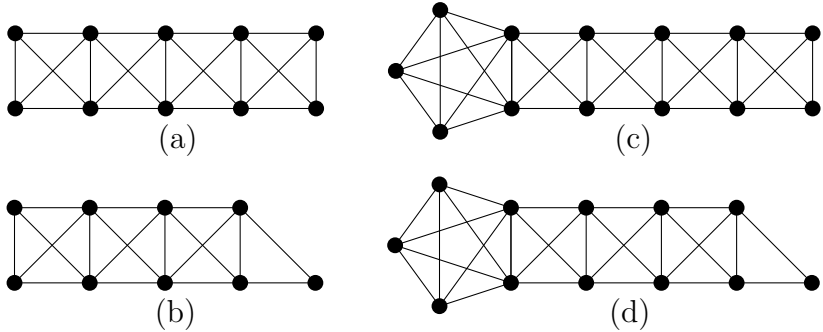
\begin{figure}[ht!]
\begin{center}
\begin{tikzpicture}
\clip(-7.28,1) rectangle (4,5.5);
\draw (-7,5)-- (-7,4);
\draw (-7,4)-- (-6,4);
\draw (-6,4)-- (-6,5);
\draw (-6,5)-- (-7,5);
\draw (-6,5)-- (-5,5);
\draw (-5,5)-- (-5,4);
\draw (-5,4)-- (-6,4);
\draw (-5,5)-- (-4,5);
\draw (-4,5)-- (-4,4);
\draw (-4,4)-- (-5,4);
\draw (-4,5)-- (-3,5);
\draw (-3,4)-- (-3,5);
\draw (-3,4)-- (-4,4);
\draw (-7.01,2.68)-- (-7.01,1.68);
\draw (-7.01,1.68)-- (-6.01,1.68);
\draw (-6.01,1.68)-- (-6.01,2.68);
\draw (-6.01,2.68)-- (-7.01,2.68);
\draw (-6.01,2.68)-- (-5.01,2.68);
\draw (-5.01,2.68)-- (-5.01,1.68);
\draw (-5.01,1.68)-- (-6.01,1.68);
\draw (-5.01,2.68)-- (-4.01,2.68);
\draw (-4.01,2.68)-- (-4.01,1.68);
\draw (-4.01,1.68)-- (-5.01,1.68);
\draw (-3.01,1.68)-- (-4.01,1.68);
\draw (-7,5)-- (-6,4);
\draw (-7,4)-- (-6,5);
\draw (-6,5)-- (-5,4);
\draw (-6,4)-- (-5,5);
\draw (-5,5)-- (-4,4);
\draw (-5,4)-- (-4,5);
\draw (-4,5)-- (-3,4);
\draw (-3,5)-- (-4,4);
\draw (-4.01,2.68)-- (-3.01,1.68);
\draw (-7.01,2.68)-- (-6.01,1.68);
\draw (-6.01,2.68)-- (-7.01,1.68);
\draw (-6.01,2.68)-- (-5.01,1.68);
\draw (-5.01,2.68)-- (-6.01,1.68);
\draw (-5.01,2.68)-- (-4.01,1.68);
\draw (-5.01,1.68)-- (-4.01,2.68);
\draw (-0.4,5)-- (-0.4,4);
\draw (-0.4,4)-- (0.6,4);
\draw (0.6,4)-- (0.6,5);
\draw (0.6,5)-- (-0.4,5);
\draw (0.6,5)-- (1.6,5);
\draw (1.6,5)-- (1.6,4);
\draw (1.6,4)-- (0.6,4);
\draw (1.6,5)-- (2.6,5);
\draw (2.6,5)-- (2.6,4);
\draw (2.6,4)-- (1.6,4);
\draw (2.6,5)-- (3.6,5);
\draw (3.6,4)-- (3.6,5);
\draw (3.6,4)-- (2.6,4);
\draw (-0.41,2.68)-- (-0.41,1.68);
\draw (-0.41,1.68)-- (0.59,1.68);
\draw (0.59,1.68)-- (0.59,2.68);
\draw (0.59,2.68)-- (-0.41,2.68);
\draw (0.59,2.68)-- (1.59,2.68);
\draw (1.59,2.68)-- (1.59,1.68);
\draw (1.59,1.68)-- (0.59,1.68);
\draw (1.59,2.68)-- (2.59,2.68);
\draw (2.59,2.68)-- (2.59,1.68);
\draw (2.59,1.68)-- (1.59,1.68);
\draw (3.59,1.68)-- (2.59,1.68);
\draw (-0.4,5)-- (0.6,4);
\draw (-0.4,4)-- (0.6,5);
\draw (0.6,5)-- (1.6,4);
\draw (0.6,4)-- (1.6,5);
\draw (1.6,5)-- (2.6,4);
\draw (1.6,4)-- (2.6,5);
\draw (2.6,5)-- (3.6,4);
\draw (3.6,5)-- (2.6,4);
\draw (2.59,2.68)-- (3.59,1.68);
\draw (-0.41,2.68)-- (0.59,1.68);
\draw (0.59,2.68)-- (-0.41,1.68);
\draw (0.59,2.68)-- (1.59,1.68);
\draw (1.59,2.68)-- (0.59,1.68);
\draw (1.59,2.68)-- (2.59,1.68);
\draw (1.59,1.68)-- (2.59,2.68);
\draw (-0.4,4)-- (-0.4,5);
\draw (-0.4,5)-- (-1.35,5.31);
\draw (-1.35,5.31)-- (-1.94,4.5);
\draw (-1.94,4.5)-- (-1.35,3.69);
\draw (-1.35,3.69)-- (-0.4,4);
\draw (-0.41,1.68)-- (-0.41,2.68);
\draw (-0.41,2.68)-- (-1.36,2.99);
\draw (-1.36,2.99)-- (-1.95,2.18);
\draw (-1.95,2.18)-- (-1.36,1.37);
\draw (-1.36,1.37)-- (-0.41,1.68);
\draw [color=black](-5.25,3.94) node[anchor=north west] {(a)};
\draw [color=black](-5.29,1.62) node[anchor=north west] {(b)};
\draw [color=black](0.82,3.94) node[anchor=north west] {(c)};
\draw [color=black](0.75,1.62) node[anchor=north west] {(d)};
\draw (-0.4,5)-- (-1.94,4.5);
\draw (-0.4,5)-- (-1.35,3.69);
\draw (-1.35,5.31)-- (-1.35,3.69);
\draw (-1.35,5.31)-- (-0.4,4);
\draw (-1.94,4.5)-- (-0.4,4);
\draw (-1.36,2.99)-- (-0.41,1.68);
\draw (-1.36,2.99)-- (-1.36,1.37);
\draw (-1.95,2.18)-- (-0.41,2.68);
\draw (-1.95,2.18)-- (-0.41,1.68);
\draw (-1.36,1.37)-- (-0.41,2.68);
\fill [color=black] (-7,5) circle (3pt);
\fill [color=black] (-7,4) circle (3pt);
\fill [color=black] (-6,4) circle (3pt);
\fill [color=black] (-6,5) circle (3pt);
\fill [color=black] (-5,5) circle (3pt);
\fill [color=black] (-5,4) circle (3pt);
\fill [color=black] (-4,5) circle (3pt);
\fill [color=black] (-4,4) circle (3pt);
\fill [color=black] (-3,5) circle (3pt);
\fill [color=black] (-3,4) circle (3pt);
\fill [color=black] (-7.01,2.68) circle (3pt);
\fill [color=black] (-7.01,1.68) circle (3pt);
\fill [color=black] (-6.01,1.68) circle (3pt);
\fill [color=black] (-6.01,2.68) circle (3pt);
\fill [color=black] (-5.01,2.68) circle (3pt);
\fill [color=black] (-5.01,1.68) circle (3pt);
\fill [color=black] (-4.01,2.68) circle (3pt);
\fill [color=black] (-4.01,1.68) circle (3pt);
\fill [color=black] (-3.01,1.68) circle (3pt);
\fill [color=black] (-0.4,5) circle (3pt);
\fill [color=black] (-0.4,4) circle (3pt);
\fill [color=black] (0.6,4) circle (3pt);
\fill [color=black] (0.6,5) circle (3pt);
\fill [color=black] (1.6,5) circle (3pt);
\fill [color=black] (1.6,4) circle (3pt);
\fill [color=black] (2.6,5) circle (3pt);
\fill [color=black] (2.6,4) circle (3pt);
\fill [color=black] (3.6,5) circle (3pt);
\fill [color=black] (3.6,4) circle (3pt);
\fill [color=black] (-0.41,2.68) circle (3pt);
\fill [color=black] (-0.41,1.68) circle (3pt);
\fill [color=black] (0.59,1.68) circle (3pt);
\fill [color=black] (0.59,2.68) circle (3pt);
\fill [color=black] (1.59,2.68) circle (3pt);
\fill [color=black] (1.59,1.68) circle (3pt);
\fill [color=black] (2.59,2.68) circle (3pt);
\fill [color=black] (2.59,1.68) circle (3pt);
\fill [color=black] (3.59,1.68) circle (3pt);
\fill [color=black] (-1.35,5.31) circle (3pt);
\fill [color=black] (-1.94,4.5) circle (3pt);
\fill [color=black] (-1.35,3.69) circle (3pt);
\fill [color=black] (-1.36,2.99) circle (3pt);
\fill [color=black] (-1.95,2.18) circle (3pt);
\fill [color=black] (-1.36,1.37) circle (3pt);
\end{tikzpicture}
\caption{Graphs $H_5$  (a), $H_5^-$ (b), $\Lambda_{5,5}$ (c), and $\Lambda_{5,5}^-$ (d)}
\label{fig2}
\end{center}
\end{figure}

\proof
Let $H_r = P_r \boxtimes K_2$, $r\ge 4$, and let $H_r^-$ be the graph obtained from $H_r$ by removing a vertex of degree $3$. Let $\Lambda_{k,r} = A(H_r,e;K_k, f)$, where $e$ is an edge of $H_r$ both of its endpoints are of degree $3$, and $f$ is an arbitrary (but fixed) edge of $K_k$. The graph $\Lambda_{k,r}^- = A(H_r^-,e;K_k, f)$ is defined analogously. See Fig.~\ref{fig2} where the graphs $H_5$,  $H_5^-$, $\Lambda_{5,5}$, and $\Lambda_{5,5}^-$ are presented. 

By Proposition~\ref{prop:amalg}, each of the graphs $H_r$, $H_r^-$, $\Lambda_{k,r}$, and
$\Lambda_{k,r}^-$ is a max-mdim graph. If $n-t+1$ is even, then the graph $\Lambda_{t-1,\frac{n-t+3}{2} }$ is of maximum degree $t$, while if  $n-t+1$ is odd, the graph $\Lambda_{t-1,\frac{n-t+4}{2} }^-$ is of maximum degree $t$. 

To complete the argument note that $n(\Lambda_{t-1,\frac{n-t+3}{2} }) = (t-1) + 2\left(\frac{n-t+3}{2}-1\right) = n$ and $n(\Lambda_{t-1,\frac{n-t+4}{2} }^-) = (t-1) + 2\left(\frac{n-t+4}{2}-1\right)-1 = n$.
\qed

At the beginning of the section we have observed that a graph obtained from a complete graph by removing a matching is a max-mdim graph provided that it contains at least two universal vertices. This fact generalizes as follows. 

\begin{prop}
If $G$ is a graph, then the following holds. 
\begin{enumerate}
\item[(i)] If $G$ has at least two universal vertices, then $G$ is a max-mdim graph.
\item[(ii)] 
If $G$ has exactly one universal vertex, then $\mdim(G)=n(G)-1$.
\end{enumerate}
\end{prop}

\proof
(i) Let $x$ and $y$ be arbitrary universal vertices of $G$. If $u\in V(G)\setminus \{x,y\}$, then $x$ (or $y$ for that matter) is a maximal neighbor of $u$. Moreover, $x$ is a maximal neighbor of $y$, and $y$ is  a maximal neighbor of $x$. By Theorem~\ref{thm:motivation}, $G$ is a max-mdim graph. 

(ii) Assume now that $x$ is the unique universal vertex of $G$. Let $W$ be a mixed resolving set for $G$. By Lemma~\ref{lem:max-neigbor} we get $\mdim(G) \geq n(G)-1$. To complete the argument we claim that $V(G)\backslash\{x\}$ is a mixed resolving set for $G$. To do this, let $\{a,b\}\subseteq\,V(G)\cup\,E(G)$. If $a\in\,V(G)\backslash\{x\}$ and $b=ax$, then $\deg_G(a) \leq n(G)-2$. Thus $V(G)\backslash N_G[a]\neq\emptyset$ and $2=d_G(a,v) \neq d_G(ax,v)=1$ for each $v\in V(G) \backslash N_G[a]$. Otherwise, there exists $u \in V(G)\backslash\{x\}$ such that $d_G(a,u)=0$ and $d_G(b,u)\geq 1$, or $d_G(a,u)\geq1$ and $d_G(b,u)=0$. Therefore, $V(G)\backslash\{x\}$ is a mixed resolving set for $G$.
\qed

\section{Graphs with cut vertices}
\label{sec:cut}
  
In this section we consider the mixed metric dimension of graphs $G$ with cut vertices and bound from the above their mixed metric dimension by $n(G) - \zeta(G)$. This of course implies (as we already know) that no graph with a cut vertex is a max-mdim graph. As a consequence we determine the  mixed metric dimension of block graphs. 

\begin{thm}\label{thm:cut}
If $W$ is a mixed resolving set of a graph $G$, then the following holds.
\begin{enumerate}
\item[(i)] If $v$ is a cut vertex of $G$, then $W\backslash\{v\}$ is a mixed resolving set of $G$. 
\item[(i)] $\mdim(G) \le  n-\zeta(G)$. Moreover, equality holds if and only if each vertex from $V(G)\setminus \CV(G)$ has a maximal neighbor in $G$.
\end{enumerate} 
\end{thm}

\proof
(i) If $v\not\in W$, then we have nothing to prove, hence assume in the remainder that $v\in W$. Let $G_1,\ldots,G_k$, $k\ge 2$, be the components of $G-v$, and for each $i\in [k]$ select a neighbor $v_i$ of $v$ in $G_i$. Since $v$ is a cut vertex, $d_G(vv_i,x)=d_G(v,x)$ for each $x\in V(G)\backslash V(G_i)$. As also $d_G(vv_i,v) = d_G(v,v) = 0$, there must be a vertex in $W\cap V(G_i)$ that distinguishes $vv_i$ and $v_i$. For each $i\in [k]$ select such a vertex $w_i$.

Consider now two arbitrary elements $a$ and $b$ from $V(G) \cap E(G)$. Assume first that $a$ and $b$ belong to some $G_i$. If $d_G(a,v) = d_G(b,v)$, then $a$ and $b$ are necessarily resolved by some vertex from $W\cap V(G_i)$. On the other hand, if $d_G(a,v) \ne d_G(b,v)$, then $a$ and $b$ are resolved by each $w_j$, where $j\ne i$. Assume next that $a$ lies in $G_i$ and $b$ belongs to $G_j$, where $i\ne j$. If $d_G(a,v) = d_G(b,v)$, then $a$ and $b$ are resolved by some vertex from $W\setminus \{v\}$. Assume that $d_G(a,v) \ne d_G(b,v)$, let  without loss of generality $d_G(a,v) > d_G(b,v)$ holds. Then we claim that $a$ and $b$ are resolved by $w_j$. Indeed, suppose on the contrary that $d_G(a,w_j) = d_G(b,w_j)$. Then
$$d_G(a,w_j) = d_G(a,v) + d_G(v,w_j) = d_G(b,w_j) \le d_G(b,v) + d_G(v,w_j)\,,$$
which in turn implies that $d_G(a,v) \le d_G(b,v)$, a contradiction. We have thus proved that each pair of elements from $V(G) \cap E(G)$ is resolved by some vertex from $W\setminus  \{v\}$, hence (i) holds. 

(ii) Since a mixed metric basis is a mixed resolving set of smallest cardinality, the inequality $\mdim(G) \le  n-\zeta(G)$  follows immediately from (i). To prove the equality part, suppose first that that each vertex from $V(G)\setminus \CV(G)$ has a maximal neighbor. Then Lemma~\ref{lem:max-neigbor} together with the already proved inequality $\mdim(G) \le  n-\zeta(G)$ yields $\mdim(G)= n-\zeta(G)$. Conversely, suppose that $\mdim(G)= n-\zeta(G)$ and suppose on the contrary that $v\in V(G)\setminus \CV(G)$ has no maximal neighbor in $G$. Then we claim that $V(G) \setminus (\CV(G)\cup \{v\})$ is a mixed resolving set for $G$. Indeed, we already know that $V(G)\setminus \CV(G)$ is a a mixed resolving set, so the only problem could be that a neighbor $u$ of $v$ would not be distinguished from the edge $uv$, because $G-v$ is a connected graph and by (i), $V(G)\setminus(\CV(G)\cup\{v\})$ is a mixed resolving set for it. However, since $v$ has no maximal neighbor, there exists $x\in N_G(v)$ such that $ux\notin E(G)$. But then $d_G(x,uv)=1$ and $d_G(x,u) = 2$. Hence $V(G) \setminus (\CV(G)\cup \{v\})$ is a mixed resolving set, a contradiction to the assumption that $\mdim(G)= n-\zeta(G)$.
\qed

Clearly, no cut vertex can have a maximal neighbor. Hence the equality part of Theorem~\ref{thm:cut}(ii) can be rephrased by saying that $\mdim(G) = n(G)$ if and only if every vertex of the graph $G$ has a maximal neighbor, which is of course Theorem~\ref{thm:motivation}. Theorem~\ref{thm:cut} also implies the following. 

\begin{cor}
\label{cor:block}
If $G$ is a block graph, then $\mdim(G) = n-\zeta(G)$.
\end{cor}

\proof
Just observe that if  $v\in V(G)\setminus \CV(G)$, then $v$ is a simplicial vertex and hence clearly has a maximal neighbor in $G$. The result then follows from Theorem~\ref{thm:cut}(ii). 
\qed

Corollary~\ref{cor:block} in turn implies that if $T$ is a tree, then $\mdim(T)$ is the number of the leaves of $T$, a result first proved in~\cite[Theorem 4.3]{KKTY-2017}.

\section*{Acknowledgments}
Sand Klav\v{z}ar was supported by the Slovenian Research Agency (ARRS) under the grants P1-0297, J1-2452, N1-0285.

\section*{Declaration of interests}

The authors declare that they have no known competing financial interests or personal relationships that could have appeared to influence the work reported in this paper.

\section*{Data availability}

Our manuscript has no associated data.


\end{document}